\documentclass[a4paper,12pt]{article}

\usepackage{amssymb,amsmath}
\usepackage{latexsym}
\usepackage{graphicx}
\newtheorem{theorem}{Theorem}[section]
\newtheorem{lem}[theorem]{Lemma}

\newtheorem{conj}[theorem]{Conjecture}

\newtheorem{claim}[theorem]{Claim}

\newcommand{\ep}{\hfill$\Box$\end{proof}\medskip}

\title{Edge-decomposition of graphs into copies of a tree with four edges}
\author{J\'anos Bar\'at \thanks{Current address: School of Mathematical
Sciences, Monash University, Victoria 3800 Australia. Research is supported by OTKA Grants PD~75837 and K~76099, and the J\'anos Bolyai Research Scholarship of the
Hungarian Academy of Sciences.}\\
\small Department of Computer Science and Systems Technology\\[-0.8ex]
\small University of Pannonia, Egyetem u.\ 10, 8200 Veszpr\'em, Hungary\\
\small \texttt{barat@dcs.vein.hu}\\
and\\
D\'aniel Gerbner\thanks{Research supported in part by the
Hungarian NSF, under contract NK 78439.}\\
\small Hungarian Academy of Sciences, Alfr\'ed R\'enyi Institute\\[-0.8ex]
\small of Mathematics, P.O.B. 127, Budapest H-1364, Hungary\\
\small \texttt{gerbner@renyi.hu}
}
\date{\today}

\begin{document}

\maketitle

\begin{abstract}
We study edge-decompositions of highly connected graphs into copies of a given tree.
In particular we attack the following conjecture by Bar\'at and Thomassen:
for each tree $T$, there exists a natural
number $k_T$ such that if $G$ is a
$k_T$-edge-connected graph, and $|E(T)|$ divides $|E(G)|$,
then $E(G)$ has a decomposition into copies of $T$.
As one of our main results it is sufficient to prove the conjecture for bipartite graphs.
Let $Y$ be the unique tree with degree sequence $(1,1,1,2,3)$.
We prove that if $G$ is a $191$-edge-connected graph of size divisible by $4$, then $G$
has a $Y$-decomposition.
This is the first instance of such a theorem, in which the tree is different from a path or a star.
\end{abstract}

MSC: 05C40, 05C70

\section{Introduction}

Our notations and concepts strictly follow \cite{thomassen2}.
A graph $G$ has an {\it $H$-de\-com\-po\-si\-tion}, if the edges of $G$
can be decomposed into subgraphs isomorphic to $H$.
There is a necessary condition: $|E(H)|$ divides $|E(G)|$.
In what follows, we always assume this hypothesis.
The general problem of $H$-decompositions
was proved to be NP-complete for any $H$ of size greater than 2 by
Dor and Tarsi \cite{dortarsi}.
However Bar\'at and Thomassen \cite{bt} posed the following

\begin{conj}\label{conj1}
For each tree $T$, there exists a natural
number $k_T$ such that the following holds: if $G$ is a
$k_T$-edge-connected graph such that $|E(T)|$ divides $|E(G)|$,
then $G$ has a $T$-decomposition.
\end{conj}

In Section~\ref{enough} we prove that it is sufficient to prove the conjecture for bipartite graphs.

\medskip

\noindent {\bf Theorem} \it
Let $T$ be a tree with $t$ edges, where $t>3$. The following two statements are equivalent.\\
(i) There exists a natural number $k_T$ such that for any $k_T$-edge-connected bipartite graph $G$,
 if $t$ divides $|E(G)|$, then $G$ has a $T$-decomposition.\\
(ii) There exists a natural number $k'_T$ such that for any $k'_T$-edge-connected graph $G$,
if $t$ divides $|E(G)|$, then $G$ has a $T$-decomposition. \rm

\medskip

In many cases $k$-edge-connectivity is provided by the existence of $k$ edge-disjoint spanning trees.
Nash-Williams \cite{n-w} and Tutte \cite{tu} independently proved the following converse.

\begin{theorem}\label{nwt}
 If $k$ is a natural number, and $G$ is a $2k$-edge-connected graph,
then $G$ contains $k$ pairwise edge-disjoint spanning trees.
\end{theorem}

At the time of posing there was no tree with at least three edges, for which Conjecture~\ref{conj1} was
known to be true.
A nice and thorough introduction to the subject is \cite{thomassen1},
where Thomassen proved that every $207$-edge-connected graph $G$ has a set $E$ of at most $6$
edges such that $G-E$ has a $4$-path-decomposition.
Approximately the same time Thomassen \cite{thomassen2} proved

\begin{theorem}\label{3path} If $G$ is a $171$-edge-connected graph of size divisible
by $3$, then $G$ has a $3$-path-decomposition.
\end{theorem}

Very recently Thomassen \cite{ct3} proved the following

\begin{theorem}
Let $k$ be any natural number. If $G$ is a $(2k^2+k)$-edge-connected graph and $|E(G)|$ is divisible by $k$,
then $G$ has a $K_{1,k}$-decomposi\-tion.
\end{theorem}

It is also mentioned in \cite{ct3} that Conjecture~\ref{conj1} holds for any path with $2^t$ edges.
In view of these results by Thomassen our Theorem~\ref{fo} is the first confirmation of Conjecture~\ref{conj1}, where
$T$ is different from a path or a star.

The proof of Theorem~\ref{3path} consists of three main ingredients.
In principle the method could be applied to any tree $T$.
Let $G$ be a graph of sufficiently high edge-connectivity,
and let $T$ be a tree on $k$ edges.
In a nutshell Thomassen set up the following scheme:

\noindent 1. Remove copies of $T$ from $G$ such that a bipartite graph $G[A,B]$ remains that still contains many edge-disjoint spanning trees.\\
2. Remove more copies of $T$ such that each degree in $A$ becomes divisible by $k$,
and the rest still contains some edge-disjoint spanning trees.\\
3. Group the edges from $A$ such that copies of $T$ arise, which altogether decompose
the rest.

It looks apparent to use the number of edge-disjoint spanning trees in a condition
instead of edge-connectivity.
We formulate a tempting qualitative conjecture, which is partly inspired by Thomass\'e \cite{stef}.

\begin{conj}\label{spt}
Let $T$ be a tree with $t$ edges.
If a graph $G$ contains $t$ edge-disjoint spanning trees and $t$ divides $|E(G)|$,
then $G$ has a $T$-decomposition.
\end{conj}

Using a maximum cut idea it is easy to prove the following

\begin{lem}\label{bipart} If $k$ is a natural number and $G$ is a
$(2k-1)$-edge-connected graph, then $G$ has a bipartition such
that $G[A,B]$ is $k$-edge-connected.
\end{lem}

%

We make heavy use of the following result by
Ellingham, Nam and Voss \cite{elli}.

\begin{lem}\label{env}
If $G$ is an $m$-edge-connected graph, then $G$ has a spanning tree $T$ such that $d_T(v) \le \lceil d_G(v)/m\rceil +2$.
\end{lem}

There is a unique tree with degree sequence $(1,1,1,2,3)$.
For simplicity we call it $Y$.
The vertex of degree 3 in $Y$ is the 3-vertex, and the vertex of degree 2 is the $2$-vertex for any further reference.
In Section~\ref{main} we prove Theorem~\ref{fo} that every $191$-edge-connected graph $G$ has a $Y$-decomposition.

\begin{figure}[ht]
 \begin{center}
  \includegraphics[scale=0.5]{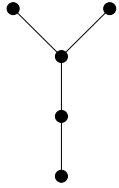}
 \end{center}
  \caption{The unique tree on 5 vertices with a vertex of degree 3}
\label{y_fa}
\end{figure}

These constants in the edge-connectivity are most likely far from optimal. 
At the end of the paper we list a few examples indicating some lower bounds.

\section{Making the graph bipartite} \label{enough}

In Thomassen's scheme the first step
is to delete some copies of the tree such that
the remaining graph is a highly edge-connected bipartite graph.
It was mentioned in \cite{thomassen2} that perhaps this method works
for every tree.
In this section we validate this hypothesis.
We need the following result that is practically a consequence of
Lemma~\ref{env}.

\begin{lem}\label{uj}
For any natural numbers $k,\ell$ and $m$, where $m>3$, if $G$ contains $km^{2\ell}$ edge-disjoint spanning trees, then we can choose subgraphs
$M_1\subset M_2\subset \dots \subset M_{\ell}\subset M_{\ell+1}$ such that $M_i$ is the union of $km^{2(i-1)}$ edge-disjoint spanning trees and $d_{M_{i}}(v) \le d_{M_{i+1}}(v)/m$
for every vertex $v$ and $1\le i\le \ell$.
\end{lem}

\begin{proof}
By Lemma~\ref{env} if we are given $m^2$ edge-disjoint spanning trees of $G$, then there exists a spanning tree $T$
such that $d_{T}(v) \le \lceil d_*(v)/m^2\rceil+2$, where $d_*(v)$ is the total degree in that particular collection of $m^2$ spanning trees.
Now $\lceil d_*(v)/m^2\rceil+2\le d_*(v)/m$ since $d_*(v)\ge m^2$ and $m>3$.

We prove the Lemma by induction on $\ell$.
We start with the base $\ell=1$.
If we are given $km^2$ edge-disjoint spanning trees, then we divide them
into $k$ equal sets $S_1,\dots,S_k$ each of them being the union of $m^2$ spanning trees.
By the previous argument there exists a spanning tree $T_i$ in $S_i$
such that $d_{T_i}(v)\le d_{S_i}(v)/m$ for any vertex $v$, where $1\le i\le k$.
Now $M_1:=\cup_{i=1}^kT_i$ and $M_2:=\cup_{i=1}^kS_i$.
Summing all the inequalities yields $d_{M_1}(v)\le d_{M_2}(v)/m$ as required for $\ell=1$.

In the induction step let $\ell>1$.
There are $m^2km^{2(\ell-1)}$ edge-disjoint spanning trees given.
We partition this set of spanning trees into $S_1,\dots,S_{km^{2(\ell-1)}}$ such that each $S_i$ contains $m^2$ spanning trees.
For every $i$ there exists a spanning tree $T_i$ in $S_i$
such that $d_{T_i}(v)\le d_{S_i}(v)/m$, where $1\le i\le k$.
By the induction hypothesis we can find $M_1\subset M_2\subset \dots \subset M_{\ell}$ in $\cup_{i=1}^{km^{2(\ell-1)}}T_i$.
Finally let $M_{\ell+1}=\cup_{i=1}^{km^{2(\ell-1)}}S_i$.
Now $d_{M_{i}}(v) \le d_{M_{i+1}}(v)/m$ is satisfied for every vertex $v$ and $1\le i\le \ell-1$
by the induction hypothesis and for $i=\ell$ by the same summation as in the base case.
\ep

In our paper we only apply the disjoint union of spanning trees in $M_i$ to ensure sufficiently large minimum degree.

Next we show that it is sufficient to prove Conjecture~\ref{conj1} for bipartite graphs:

\begin{theorem} \label{paros} Let $T$ be a tree with $t$ edges, where $t>3$. The following two statements are equivalent.\\
(i) There exists a natural number $k_T$ such that for any $k_T$-edge-connected bipartite graph $G$,
 if $t$ divides $|E(G)|$, then $G$ has a $T$-decomposition.\\
(ii) There exists a natural number $k'_T$ such that for any $k'_T$-edge-connected graph $G$,
if $t$ divides $|E(G)|$, then $G$ has a $T$-decomposition.
\end{theorem}

\begin{proof} We only prove the non-trivial implication.
Let $k'_T=8t^{2t+3}+4k_T-1$.\linebreak
By Lemma~\ref{bipart} we can find a partition $(A,B)$ of the
vertex set such that the edge-connectivity of $G[A,B]$ is $4t^{2t+3}+2k_T$.
By Theorem~\ref{nwt}
there are at least $2t^{2t+3}+k_T$ pairwise edge-disjoint spanning trees in $G[A,B]$.
In what follows we show how to delete all edges inside $A$ by removing copies of $T$ from $G[A,B]$ using
at most $t^{2t+3}$ of the spanning trees. We can repeat the same procedure to empty $B$.
After that the remaining $k_T$ spanning trees provide $k_T$-edge-connectivity of the remaining bipartite graph.

First we arbitrarily delete copies of $T$ from $G[A]$ as long as possible.
We partition the remaining edges into subgraphs of $T$ as follows.
In each round we find a subtree of $T$ and remove it from the current graph.
After finishing a round we start another one on the remaining edges of $G[A]$.
Let us make a breadth-first search in $T$ starting at a vertex $x_1$ at level 0.
Select a non-isolated vertex $v_1$ of $A$ to play the role of vertex $x_1$ of $T$.
We try to fold $T$ into the remainder of $G[A]$ by the levels of the breadth-first search.
We first consider the candidates for the vertices in depth 1.
If the degree of $x_1$ is $d$ in $T$, then select $d$ edges incident to $v_1$ to form a depth 1 copy of $T$.
If the degree of $x_1$ is too large, then select all edges incident to $v_1$ to form a partial copy.
Next we consider the neighbours of $v_1$ trying to attach vertices to the current leaves to create a
depth 2 (partial) copy of $T$.
We continue this process level by level.
We have to skip those edges of $G[A]$ that would create a cycle in the (partial) copy of $T$.
At some point we run out of possible extensions since there are no complete copies of $T$ in the present graph.
In that case we finish the round with the subgraph that we started at $v_1$.
We repeat this process on the set of remaining edges starting at another arbitrary non-isolated vertex.
In this way we create a set $\mathcal{H}$ of subgraphs of $G[A]$, each of which is isomorphic to
a subtree of $T$.

Let $H$ be a graph in $\mathcal{H}$.
A vertex $v$ is \emph{unsaturated} in $H$, if it corresponds to a
vertex of degree $d$ in $T$, but the degree of $v$ in $H$ is
smaller than $d$.
For every vertex $v$ in $A$ there are at
most $t$ trees in $\mathcal{H}$, in which $v$ is unsaturated. Indeed
consider the first occasion $v$ has too few available edges to copy the next level of $T$ in the above process.
At the end of the round we get a subgraph $H$, in which $v$ is unsaturated.
At this point all remaining edges incident to $v$ must be incident to other vertices in $H$.
We could not use them, because any of them created a cycle.
The size of $H$ is at most $t-1$, and later we probably used these edges at $v$ one by one to create partial copies
of $T$.
Therefore, $v$ can be unsaturated in at most $t-1$ members of $\mathcal{H}$.
With a different argument we can prove this to be at most $t/2$, but for simplicity we use $t-1$ in the counting below.

\begin{claim}
The set $\mathcal{H}$ can be partitioned into $t^2$ sets
$\mathcal{H}_1, \dots, \mathcal{H}_{t^2}$ such that for each $i$ and each
vertex $v$ in $A$, there is at most one tree in $\mathcal{H}_i$  where
$v$ is unsaturated, $1\le i\le t^2$.
\end{claim}

\begin{proof}
As long as possible select trees into $\mathcal{H}_1$ without violating the property.
In the rest create similarly $\mathcal{H}_2, \dots ,\mathcal{H}_{t^2}$.
If a tree $H$ is not selected in $\mathcal{H}_i$, then one of its vertices violates the property, say $v$.
Necessarily $v$ is unsaturated in a tree already selected to be in $\mathcal{H}_i$.
The size of $H$ is at most $(t+1)$ and $v$ is unsaturated at most $(t-1)$ times.
Therefore, such a bad event can happen at most $(t+1)(t-1)$ times and the claim holds.
\ep

Back to the proof of the Theorem:
We save $k_t$ edge-disjoint spanning trees for later, and
partition the remaining available spanning trees into $2t^2$ sets each of them containing $t^{2t+1}$ spanning trees.
Consider $\mathcal{H}_1$ and $t^{2t+1}$ spanning trees.
We apply Lemma~\ref{uj} with $t=k=\ell=m$ to the graph being the union of these spanning trees.
Here $\ell=t-1$ would suffice.
We get the subgraphs $M_1,\dots, M_{t+1}$ with the properties in Lemma~\ref{uj}.

We first describe our process for one
incomplete copy $H\in \mathcal{H}_1$ of the tree $T$.
Select all unsaturated vertices of $H$ forming a set $D_0(H)$.
We extend $H$ into a copy of $T$ level by level starting at $D_0(H)$.
In the first step we use the edges from the $t$ edge-disjoint spanning trees that
form $M_1$.
Notice that every vertex in $D_0(H)$ has degree at least $t$ in $M_1$.
Since at most $t$ edges are to be added in total, we can choose to extend the
tree $H$ to $H(1)$ by using edges with pairwise distinct endvertices in $B$.
These new vertices form the set $D_1(H)$ and we produced a larger subtree $H(1)$ of $T$.
We next use the edges of $M_2$ to proceed with the completion of $T$.
We know that $d_{M_2}(v)\ge td_{M_1}(v)$.
Therefore, we can choose the edges to complete the next level of $T$ such that
their endvertices (at most $t$) are pairwise disjoint, and
they are distinct from the (at most $t$) vertices in $H(1)$.
Generally in step $(j+1)$ we can use the edges of $M_{j+1}\setminus M_{j}$ to extend
$H(j)$ to $H(j+1)$. Since $d_{M_{j+1}}(v)\ge td_{M_{j}}(v)$, we can avoid creating a cycle
and we can fulfil any demand for an edge.

Instead of a concrete $H$ we do the above process simultaneously for all subtrees in $\mathcal{H}_1$.
Let all unsaturated vertices in $\mathcal{H}_1$ form $D_0(\mathcal{H}_1)$.
We extend each subtree to a copy of $T$.
Each vertex of $A$ is unsaturated at most once in $\mathcal{H}_1$ and $d_{M_1}(v)\ge t$ for any vertex $v$.
Therefore, $M_1$ provides sufficient amount of edges to add new vertices forming $D_1(\mathcal{H}_1)$
and to create the set of subtrees $\mathcal{H}_1(1)$ to complete level 1.
It might happen that many edges from different vertices of $A$ connect to a particular vertex $u$ of $B$.
It makes $d_{M_1}(u)$ large.
Anyway the growth $d_{M_2}(u)\ge td_{M_1}(u)$ ensures us that we can continue with step 2.

In general in step $(j+1)$ we extend the trees $H(j)$ to $H(j+1)$.
We form $D_{j+1}(\mathcal{H}_1)$
and create $\mathcal{H}_1(j+1)$ to complete level $(j+1)$.
Consider a vertex $v$, which is on level $j$ in $H(j)$ and possibly in some other trees.
We know that $v$ is on level $j$ in at most $d_{M_j}(v)$ trees. 
We extend the notion of unsaturated vertices: a vertex $v$ of $H(j)\setminus H(j-1)$ is {\it unsaturated} 
if its degree in $H(j)$ (which is exactly 1) is smaller than the degree of the corresponding vertex in $T$. 
For each subtree where $v$ is unsaturated we select $(t-1)$ edges of $M_{j+1}\setminus M_j$ incident to $v$.
As $d_{M_{j+1}}(v)\ge td_{M_j}(v)$ we can make the selection such that
every edge is chosen at most once. 
We do the same for every vertex in $D_j(\mathcal{H}_1)$.

We use only the selected edges for the extension of the subtrees. 
We have to show that any subtree $H(j)$ in $\mathcal{H}_1(j)$ can be extended to $H_1(j+1)$ using those edges. 
Now there might be many unsaturated vertices $v_1, \dots, v_q$ in $H(j)$, but at most $(t-1)$ edges have to be added altogether.
As we selected $t-1$ edges for every $v_i$, we can choose the edges of the next level such that no cycle occurs. 
Two different subtrees $H(j)$ and $H'(j)$ can be extended simultaneously since the selected sets of edges are disjoint.
Therefore, all subtrees of $\mathcal{H}_1(j)$ can be extended to $\mathcal{H}_1(j+1)$ simultaneously.

The proof is completed by repeating the argument for each
$\mathcal{H}_i$ with another sequence $M_1\subset\dots\subset M_{t+1}$ provided by the corresponding set
of $t^{2t+1}$ edge-disjoint spanning trees.

We repeat the argument with $B$ and remove all complete copies of $T$.
It yields a bipartite graph that contains the remaining set of $k_t$ edge-disjoint spanning trees.
\ep

For any fixed tree, the above edge-connectivity condition can be largely reduced.
For any such improvement, we use the same principal argument, but we can decrease the necessary number of spanning trees,
by using the structure of the fixed tree.
In particular for the graph $Y$ we show the following.

\begin{lem} \label{jo}
If $G$ is a $(4k+23)$-edge-connected graph, then we can remove some $Y$-copies
such that a bipartite graph with
$k$ edge-disjoint spanning trees remains.
\end{lem}

\begin{proof}
By Lemma~\ref{bipart} we can find a bipartition $G[A,B]$ of $G$ that is $(2k+12)$-edge-connected.
By Theorem~\ref{nwt} we find $(k+6)$ edge-disjoint spanning trees in $G[A,B]$.
Let $T_1, T_2, T_3$ be three of them.
We remove $Y$-copies from $G[A]$ arbitrarily as long as we can.
What remains in $G[A]$ is a collection of paths, cycles, stars and subgraphs of $K_4$.
We cut each path and each cycle into paths with three edges and a possible shorter path.
We select one of the middle vertices of such a 3-path to be the 3-vertex of a future $Y$-copy.
The idea is to extend these 3-paths into $Y$-copies using $T_1$, and remove them from $G$.
For a 2-path we select one endvertex to be the candidate 3-vertex of $Y$.
For a single edge we select one endvertex to be the candidate 3-vertex of $Y$,
and the other endvertex to be the 2-vertex of $Y$.
We cut the stars into 3-stars and a remaining part, which is a 2-path or a single edge as above.
For a 3-star we select a leaf to be the 2-vertex of $Y$.
Until now any vertex in $A$ is selected at most once.
For any subgraph $H$ of $K_4$ that is different from the previous ones, we do as follows.
We cut $H$ into paths of length at most three such that
after the above selection of 3- and 2-vertices of $Y$, each vertex is used at most once.
This is always possible with one exception, the triangle.

If a vertex of $A$ is selected to be a 3- or 2-vertex of $Y$, then
we extend the subgraph with edges of $T_1$ and $T_2$ to achieve a $Y$-copy
that we remove.
It works fine except for a single edge or a triangle.
In case of a single edge we have to add three additional edges to get a $Y$-copy.
For the vertex selected to be the 3-vertex we use edges from $T_1$ and $T_2$.
Now there exists an edge in $T_1$, $T_2$ or $T_3$ from the other end of the single edge
that avoids creating a cycle, hence it completes to a $Y$-copy that we remove.

In case of a triangle we cut it into a single edge and a 2-path.
We do as above for the single edge, and let $v$ be the vertex selected to be the 2-vertex.
For the 2-path we select $v$ to be the 3-vertex of $Y$.
Since we used one of $T_1-T_3$ for the single edge, there are two edges left
to use.
We create a $Y$-copy and remove it.

We have to execute the same process for $G[B]$, where we use three more spanning trees.
After all a bipartite graph remains that has at least $k$ edge-disjoint spanning trees.
\ep

Even if there are only $(k+5)$ spanning trees in $G[A,B]$, we can delete $Y$-copies using $5$ spanning trees
such that a $k$-edge-connected bipartite graph remains.
It requires a more detailed argument,
and implies an improvement by $4$ in the statement of the Lemma and subsequently in Theorem~\ref{fo}.

\section{Proof of the main theorem}\label{main}

We recall an implicit result from \cite{thomassen2}.
The second and third paragraph on page 291 describe a $P_4$-decomposition of a special
graph. We realised that the vertices of $A$ are used in the decomposition
in a balanced way.

\begin{lem}\label{vege}
Let $G$ be a $2$-edge-connected bipartite graph with classes $A$ and $B$.
If the degree of each vertex in $A$ is divisible by $3$, then $G$ can be decomposed
into paths with $3$ edges such that each vertex $v$ of $A$ is the
endvertex of $d(v)/3$ paths and middle vertex of $d(v)/3$ paths.
\end{lem}

We use this lemma in the finishing stage of the next
result that gives a sufficient edge-connectivity condition for $Y$-de\-com\-po\-si\-tions.

\begin{theorem} \label{fo}
Let $Y$ denote the tree with degree sequence $(1,1,1,2,3)$. If $G$
is a $191$-edge-connected graph of size divisible by $4$, then $G$
has a $Y$-decomposition.
\end{theorem}

\begin{proof}
We first apply Lemma~\ref{jo} with $k=42$. As a result we are
given a bipartite graph $G[A,B]$ with 42 edge-disjoint spanning
trees $T_1, \dots, T_{42}$.

In the next step
we delete some copies of $Y$ to make all degrees in $A$ divisible by 4.
In the first phase we achieve that all degrees are even.
Therefore, vertices in $A$ of odd degree are {\it bad}.
Let $M(1)$ be a subgraph of $G$ that is the union of 7
edge-disjoint spanning trees $T_1, \dots, T_7$. By
Lemma~\ref{env}, $m=7$, $M(1)$ has a spanning tree $T(1)$ such
that for each vertex $v$, $d_{T(1)}(v)\le \lceil d_{M(1)}(v)/7\rceil+2\le
d_{M(1)}(v)/2$, since $d_{M(1)}(v)\ge 7$. Similarly the union
$M(2)$ of 7 spanning trees $T_8, \dots, T_{14}$ contains a
spanning tree $T(2)$ such that
$d_{T(2)}(v)\le d_{M(2)}(v)/2$ for each vertex $v$. The union of $T(1)$ and $T(2)$
contains a spanning Eulerian subgraph $E_1$.

We start a walk on $E_1$ at a bad vertex $u_1$. We construct and delete
$Y$-copies as follows. Let $e_1$ be the edge
adjacent to $u_1$ in $E_1$, and let $e_2, e_3,\dots$ be the edges of $E_1$ in order.
Walking along $e_1$ and $e_2$ we are back in $A$ in a vertex $u_2$.
We continue this way till we arrive to another bad vertex $u_r$.
We selected an edge incident to $u_1$ that we later will remove.
Therefore, when $E_1$ possibly arrives to $u_1$ next time, it is no longer considered a bad vertex.
That is, $u_r\neq u_1$.
For every $i$ we consider $e_{2i-1}$, $e_{2i}$ and
two edges in $M(1) \cup M(2) \setminus E_1$ that are incident to $u_{i+1}$, where $1\le i \le r-1$.
These four edges form a copy of $Y$ that we delete.
In this way we delete an odd number of edges incident to $u_1$ and $u_r$, and an even
number of edges incident to any other vertex in $A$.
Therefore, the number of bad vertices decreases.
A vertex can appear multiple times in the above sequence, but that does not change the
parity of the degree.

Now we continue the walk along $E_1$ and do nothing until we find another pair of bad vertices.
We repeat the above process of removing $Y$-copies between the bad vertices.
Iterating these two steps we finish the Eulerian trail, and all degrees are now even.
There is a small remark that we have to make: there are enough edges in  $M(1) \cup M(2)
\setminus E_1$ to use.
Indeed, whenever the walk arrives to a vertex $v$, it means
there are two incident edges in $E_1$.
Hence we can find two more edges, as the degree of a vertex $v$ in $E_1$ is at most half of the
degree of $v$ in $M(1) \cup M(2)$.

In the second phase all degrees in $A$ are even.
Our goal is to remove some $Y$-copies to make all degrees divisible by 4.
Therefore, vertices in $A$ of degree $2\mod 4$ are considered {\it bad}.
As in the first phase we need an Eulerian spanning subgraph for our purposes.
Let $M(3)$ be a subgraph of $G$ that is the union of 9
edge-disjoint spanning trees $T_{15}, \dots, T_{23}$. By
Lemma~\ref{env}, $m=9$, $M(3)$ has a spanning tree $T(3)$ such
that for each vertex $v$, $d_{T(3)}(v)\le \lceil d_{M(3)}(v)/9\rceil+2\le
d_{M(3)}(v)/2-1$, since $d_{M(3)}(v)\ge 9$. Similarly the union
$M(4)$ of the spanning trees $T_{24}, \dots, T_{32}$ contains a
spanning tree $T(4)$ such that for each vertex $v$,
$d_{T(4)}(v)\le d_{M(4)}(v)/2-1$. The union of $T(3)$ and $T(4)$
contains a spanning Eulerian subgraph $E_2$.

On the Eulerian trail we mark the bad vertices.
We start the marking at a bad vertex $b_1$, and mark the bad vertices at the first appearance only.
We get a list $b_1,\dots ,b_r$ of bad vertices, and this list
reflects their order of first appearance on $E_2$.
This direction on $E_2$ is fixed from now on.

In what follows we remove $Y$-copies to achieve that all degrees in $A$ are divisible by $4$.
If $v$ is a bad vertex, then we remove 2 or 6 edges incident to $v$ during the process,
when we arrive to the marked copy of $v$.
If $x$ is an unmarked vertex, then we remove precisely 4 edges.
If $x$ is a vertex on $E_2$, then let $x^+$ be the next vertex of $A$ on $E_2$.
There are two building bricks:\\
1. \it remove a $Y$-copy at $x$ \rm is a step, when two consecutive edges of $E_2$ starting at $x$,
and two edges of $M(3) \cup M(4) \setminus E_2$ at $x^+$ are removed.\\
2. \it remove a reversed $Y$-copy at $x$ \rm is a step, when two consecutive edges of $E_2$ starting at $x$,
and two edges of $M(3) \cup M(4) \setminus E_2$ at $x$ are removed.

We start at $b_1$ and remove a $Y$-copy.
We continue along $E_2$ and remove all edges of $E_2$ two by two.
Every such pair of edges corresponds to a $2$-path in a $Y$-copy, where one end is the $3$-vertex.
The only decision to make is the placement of the other two edges from $M(3) \cup M(4) \setminus E_2$.
Either at the current vertex $x$ or at the subsequent vertex $x^+$.
This is actually automatic, according to the degree condition:
we either deleted 1 or 3 edges at $x$ due to the previous $Y$-copy,
and our goal might be to remove 2, 4 or 6 edges in total.
If we need to remove one more edge at $x$, then we remove a $Y$-copy.
If we need to remove three more edges at $x$, then we remove a reversed $Y$-copy.
Notice here that finishing the Eulerian trail we get back to $b_1$.
The last condition automatically removes one more edge at $b_1$, since
the remaining number of edges has to be divisible by 4.

After this process bad vertices become good, and the degrees
of good vertices are still divisible by $4$.
Here we also remark that there are enough edges in $M(3) \cup M(4)
\setminus E_2$ to use every time the walk arrives to a vertex $v$.
This again follows from the upper bound on $d_{T(3)}$ and $d_{T(4)}$.
Whenever we arrive to $v$, it means there are two edges incident to $v$ in $E_2$, and we
need two edges (or four, at most once) in $M(3) \cup M(4)
\setminus E_2$.
Therefore, we need the degree of $v$ in $M(3)
\cup M(4) \setminus E_2$ to be at least $d_{E_2}(v)+2$, which is satisfied.

We are left with a bipartite graph $M[A,B]$, where all degrees in
$A$ are divisible by 4. Let $M(5)$ be the union of 5 spanning
trees $T_{33}, \dots, T_{37}$. By Lemma~\ref{env}, $m=5$,
$M(5)$ contains a spanning tree $T(5)$ such that for each vertex
$v$, $d_{T(5)}(v)\le \lceil d_{M(5)}(v)/5\rceil+2\le
3d_{M(5)}(v)/4$, since $d_{M(5)}(v)\ge 5$.
We similarly define $M(6)$ and find $T(6)$.
Now for every vertex $v$ in $A$, the following holds:
$d_{T(5)}(v)+d_{T(6)}(v) \le 3d_M(v)/4$. For every vertex $v$ in
$A$ we put aside $1/4$ of the edges such that $T(5)$ and
$T(6)$ remains in the graph. The remaining graph $M'$ satisfies
the conditions of Lemma~\ref{vege}.

Therefore, we can decompose $M'$ into paths of length 3 such that
for a vertex $v$ with degree $4d$ in $M$ (hence degree $3d$ in the
smaller graph $M'$), there are $d$ paths starting from $v$, and $d$
paths, where $v$ is a middle vertex.
For every vertex $v$ we glue
the $d$ edges, which we put aside in the beginning of the third
phase, one by one to the $d$ paths, where $v$ is a middle vertex.
This gives us a $Y$-decomposition.
\ep

\section*{Discussion}

The edge-connectivity constants in the solved cases of Conjecture~\ref{conj1}
are seemingly far from best possible.
There is very little known about lower bounds.
For trees with three edges:
if $T$ is the $3$-path, then there is a $2$-edge-connected graph without a $3$-path-decomposition \cite{jrp}.
In \cite{bt}, there is a $4$-edge-connected graph without a $3$-star-decomposition.
In the following picture, we give a $3$-edge-connected bipartite graph with 27 edges and
without a $3$-star-decomposition.

\begin{figure}[ht]
 \begin{center}
  \includegraphics[scale=0.5]{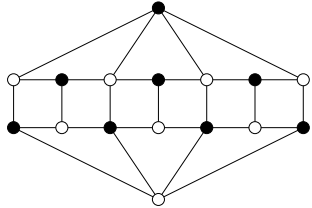}
 \end{center}
  \caption{A bipartite graph without $3$-star-decomposition}
\label{letra}
\end{figure}

There are three different trees with four edges: the $4$-star, the $4$-path and $Y$.
Consider four copies of $K_6$: the graphs $G_1,\dots, G_4$.
Add three edges between $G_i$ and $G_{i+1}$ such that we get a 6-regular graph.
This is a $6$-edge-connected graph without a $4$-star-decomposition.

If $T$ is the $4$-path, then we have the following $3$-edge-connected example
without $P_5$-decomposition.

\begin{figure}[ht]
 \begin{center}
  \includegraphics[scale=0.5]{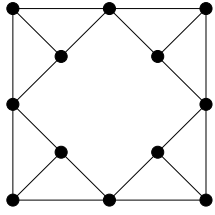}
 \end{center}
  \caption{A graph without $4$-path-decomposition}
\label{no_p5}
\end{figure}

The $4$-wheel is a $3$-edge-connected graph without a $Y$-decomposition.

One might feel that using the edge-connectivity instead of the number of spanning trees and applying
Theorem~\ref{nwt} and Lemma~\ref{bipart} is too generous.
About the sharpness of Theorem~\ref{nwt}, see \cite{cat}.

\section*{Acknowledgements}

We are greatly indebted to the organisers of  the  $2^{nd}$ Eml\'ekt\'abla Workshop January 24-27. 2011,
 Gy\"on\-gy\"os\-tarj\'an, Hungary, where the presented research  initiated.

\end{document}